\documentclass{amsart}

\usepackage{fullpage}
\usepackage{amsmath}

\usepackage{amsfonts}
\usepackage{amssymb}

\input xy
\xyoption{matrix}\xyoption{arrow}\xyoption{curve}

\begin{document}


\newcommand{\Hom}{\mathrm{Hom}}
\newcommand{\RHom}{\mathrm{RHom}^*}
\newcommand{\HOM}{\mathrm{HOM}}
\newcommand{\stHom}{\underline{\mathrm{Hom}}}
\newcommand{\Ext}{\mathrm{Ext}}
\newcommand{\Tor}{\mathrm{Tor}}
\newcommand{\HH}{\mathrm{HH}}
\newcommand{\Endo}{\mathrm{End}}
\newcommand{\ENDO}{\mathrm{END}}
\newcommand{\stEndo}{\mathrm{\underline{End}}}
\newcommand{\Tr}{\mathrm{Tr}}


\newcommand{\coker}{\mathrm{coker}}
\newcommand{\aut}{\mathrm{Aut}}
\newcommand{\op}{\mathrm{op}}
\newcommand{\add}{\mathrm{add}}
\newcommand{\ADD}{\mathrm{ADD}}
\newcommand{\ind}{\mathrm{ind}}
\newcommand{\rad}{\mathrm{rad}}
\newcommand{\soc}{\mathrm{soc}}
\newcommand{\ann}{\mathrm{ann}}
\newcommand{\im}{\mathrm{im}}
\newcommand{\chr}{\mathrm{char}}
\newcommand{\pdim}{\mathrm{p.dim}}
\newcommand{\gldim}{\mathrm{gl.dim}}


\newcommand{\rmod}{\mbox{mod-}}
\newcommand{\Rmod}{\mbox{Mod-}}
\newcommand{\lmod}{\mbox{-mod}}
\newcommand{\lMod}{\mbox{-Mod}}
\newcommand{\stmod}{\mbox{\underline{mod}-}}
\newcommand{\stlmod}{\mbox{-\underline{mod}}}

\newcommand{\gmod}[1]{\mbox{mod}_{#1}\mbox{-}}
\newcommand{\gMod}[1]{\mbox{Mod}_{#1}\mbox{-}}
\newcommand{\Bimod}[1]{\mathrm{Bimod}_{#1}\mbox{-}}

\newcommand{\proj}{\mbox{proj-}}
\newcommand{\lproj}{\mbox{-proj}}
\newcommand{\Proj}{\mbox{Proj-}}
\newcommand{\inj}{\mbox{inj-}}
\newcommand{\coh}{\mbox{coh-}}
\newcommand{\CM}{\mbox{CM}}
\newcommand{\sCM}{\mbox{\underline{CM}}}


\newcommand{\und}[1]{\underline{#1}}
\newcommand{\gen}[1]{\langle #1 \rangle}
\newcommand{\floor}[1]{\lfloor #1 \rfloor}
\newcommand{\ceil}[1]{\lceil #1 \rceil}
\newcommand{\bnc}[2]{\left(\scriptsize \begin{array}{c} #1 \\ #2 \end{array} \right)}
\newcommand{\bimo}[1]{{}_{#1}#1_{#1}}
\newcommand{\ses}[5]{\ensuremath{0 \rightarrow #1 \stackrel{#4}{\longrightarrow} 
#2 \stackrel{#5}{\longrightarrow} #3 \rightarrow 0}}
\newcommand{\A}{\mathcal{A}}
\newcommand{\B}{\mathcal{B}}
\newcommand{\uB}{\underline{\mathcal{B}}}
\newcommand{\C}{\mathcal{C}}
\newcommand{\uC}{\underline{\mathcal{C}}}
\newcommand{\D}{\mathcal{D}}
\newcommand{\E}{\mathcal{E}}
\newcommand{\F}{\mathcal{F}}
\newcommand{\p}{\mathcal{P}}
\newcommand{\ul}[1]{\underline{#1}}


\newtheorem{therm}{Theorem}[section]
\newtheorem{defin}[therm]{Definition}
\newtheorem{propos}[therm]{Proposition}
\newtheorem{lemma}[therm]{Lemma}
\newtheorem{coro}[therm]{Corollary}

\title{Periodicity of $d$-cluster-tilted algebras}
\author{Alex Dugas}
\address{Department of Mathematics, University of the Pacific, 3601 Pacific Ave, Stockton, CA 95211, USA}
\email{adugas@pacific.edu}


\begin{abstract} It is well-known that any maximal Cohen-Macaulay module over a hypersurface has a periodic free resolution of period $2$.  Auslander, Reiten \cite{DTrPer} and Buchweitz \cite{Buch} have used this periodicity to explain the existence of periodic projective resolutions over certain finite-dimensional algebras which arise as stable endomorphism rings of Cohen-Macaulay modules.  These algebras are in fact periodic, meaning that they have periodic projective resolutions as bimodules and thus periodic Hochschild cohomology as well.  The goal of this article is to generalize this construction of periodic algebras to the context of Iyama's higher AR-theory.  We let $\C$ be a maximal $(d-1)$-orthogonal subcategory of an exact Frobenius category $\B$, and start by studying the projective resolutions of finitely presented functors on the stable category $\uC$, over both $\uC$ and $\C$.  Under the assumption that $\uC$ is fixed by $\Omega^{d}$, we show that $\Omega^{d}$ induces the $(2+d)^{th}$ syzygy on $\rmod \uC$.  If $\C$ has finite type, i.e., if $\C = \add(T)$ for a $d$-cluster tilting object $T$, then we show that the stable endomorphism ring of $T$ has a quasi-periodic resolution over its enveloping algebra.  Moreover, this resolution will be periodic if some power of $\Omega^{d}$ is isomorphic to the identity on $\uC$.  It follows, in particular, that $2$-C.Y.-tilted algebras arising as stable endomorphism rings of Cohen-Macaulay modules over curve singularities, as in the work of Burban, Iyama, Keller and Reiten \cite{BIKR}, have periodic bimodule resolutions of period $4$.
\end{abstract}

\maketitle
 
\section{Introduction} 
\setcounter{equation}{0}

  In this article we describe a new way of constructing finite-dimensional endomorphism algebras with periodic Hochschild (co)homology.  In fact, we show that the endomorphism rings we consider are {\it periodic} in the sense that they have periodic projective resolutions over their enveloping algebras; i.e., $\Omega_{A^e}^n(A) \cong A$ as bimodules for some $n >0$.  Among the most notable examples of finite-dimensional algebras with this property are the preprojective algebras of Dynkin graphs, which all have period $6$.  This interesting fact was first proved by Ringel and Schofield through a calculation of the minimal projective bimodule resolutions of such algebras.  Later, Auslander and Reiten \cite{DTrPer} gave an elegant functorial argument for this periodicity, making use of the fact that these preprojective algebras can be realized as stable endomorphism rings of Cohen-Macaulay modules (in fact, as stable Auslander algebras) over $2$-dimensional simple hypersurface singularities.  Actually, their arguments establish a slightly weaker version of this periodicity, showing only that the sixth power of the syzygy functor is the identity.  Motivated by these results, Buchweitz \cite{Buch} develops the functor category arguments of Auslander and Reiten to deduce the (full) periodicity of the preprojective algebras of Dynkin graphs from the isomorphisms $\Omega^2 \cong Id$ in the corresponding stable categories of CM-modules.  More generally, his work shows how periodic algebras can arise as stable Auslander algebras of finite-type categories, and in particular as stable  endomorphism rings of $\Omega$-periodic modules.
 
  Iyama has recently developed higher-dimensional analogues of much of the classical Auslander-Reiten theory, including a theory of higher Auslander algebras \cite{Iyama1, Iyama2}.  Thus it is natural to look for generalizations of Auslander, Reiten and Buchweitz's work on periodicity to this setting.  One clue is already provided by recent work of Burban, Iyama, Keller and Reiten \cite{BIKR}, showing that symmetric algebras with $\tau$-period $2$ can be obtained as endomorphism rings of certain Cohen-Macaulay modules over $1$-dimensional hypersurface singularities.  Among the algebras they realize in this way are several algebras of quaternion type, which Erdmann and Skowronski have shown are periodic of period $4$ \cite{ErdSko}.  As Erdmann and Skowronski's result is obtained by computing minimal projective resolutions over enveloping algebras, our motivation is parallel to Buchweitz's in \cite{Buch}.  That is, we aim to generalize Buchweitz's results to explain how the $2$-periodicity of the syzygy functor in the category of CM-modules implies the $4$-periodicity of the bimodule resolutions for the appropriate endomorphism rings.
  
  It turns out that we can obtain periodic algebras more generally as endomorphism rings of periodic $d$-cluster-tilting objects in a triangulated category.  These $d$-cluster-tilting objects are in fact the objects $T$ for which $\add(T)$ satisfies Iyama's definition of a maximal $(d-1)$-orthogonal subcategory.  Hence our results are indeed analogues of Buchweitz's for Iyama's higher Auslander-Reiten theory.  We summarize our main results (see Corollary 3.1 and Theorem 3.2) in the theorem below, where $\B$ denotes an exact Frobenius category with a Hom-finite stable category $\uB$.
  
  \begin{therm} Let $T$ be a $d$-cluster tilting object in $\B$ (with $d \geq 1$) such that $\Omega^{d}T \cong T$ in $\uB$, and set $\Lambda = \Endo_{\B}(T)$ and $\Gamma = \stEndo_{\B}(T)$.  If $\Gamma$ has no semisimple blocks, then
  \begin{enumerate}
  \item $\Tor_{i}^{\Lambda}(-,\Gamma) = 0$ on $\rmod \Gamma$ for all $i \neq 0, d+1$.
  \item $\Omega^{d+2}_{\Gamma^e}(\Gamma) \cong \Tor_{d+1}^{\Lambda}(\Gamma,\Gamma) \cong \uB(T,\Omega^{d}T)$ is an invertible $(\Gamma,\Gamma)$-bimodule.  Hence $\Gamma$ has a quasi-periodic projective resolution over its enveloping algebra $\Gamma^e$.
  \item If $\Omega^{d}$ has order $r$ as a functor on $\ul{\add(T)}$, then $\Gamma$ is periodic with period dividing $(d+2)r$. 
  \end{enumerate}  
 \end{therm}
  
For $d=1$, the same conclusions were obtained by Buchweitz \cite{Buch} under the assumption (needed for (2) and (3)) that $\Lambda$ has Hochschild dimension $d+1=2$.  He then applies it to an additive generator $T$ of the finite-type category $\B = \CM(R)$ for a simple hypersurface singularity $R$ of dimension $2$ in order to deduce the periodicity of the preprojective algebras of Dynkin type.  For $d=2$, we can again take $\B = \CM(R)$ for an odd-dimensional isolated Gorenstein hypersurface (see \cite{Yosh} for instance).  Since Eisenbud's matrix factorization theorem \cite{Eis} implies that $\Omega^2 \cong Id$ on $\uB$ in this case, any $2$-cluster-tilting object in $\uB$ is automatically $2$-periodic and thus has a stable endomorphism algebra which is periodic of period $4$.  Existence of $2$-cluster-tilting objects in this setting has been studied by Burban, Iyama, Keller and Reiten \cite{BIKR}.  We will discuss this and other potential applications further in the final section.

  We typically work with right modules, unless noted otherwise.  In this case morphisms are written on the left and composed from right to left.  We also follow this convention for morphisms in abstract categories.  For a category $\A$, we shall write $\A(X,Y)$ for the set of morphisms from $X$ to $Y$ in $\A$, and we shall write $\Hom_{\A}(-,-)$ for the morphism sets in categories of functors on $\A$, such as $\rmod \A$.  Likewise $\Tor^{\A}(-,-)$ and $\Ext_{\A}(-,-)$ will be reserved for $\A$-modules.  We also follow the convention of writing $\Omega_{\A} M$ for the syzygy of an $\A$-module $M$ in order to distinguish it from the syzygy operator on $\A$ (provided this makes sense), which we write simply as $\Omega$.

\section{Functors on maximal orthogonal subcategories}
\setcounter{equation}{0}

Throughout this article, we let $k$ be a field and assume that $\B$ is an exact Krull-Schmidt, Frobenius $k$-category, which arises as a full, extension-closed subcategory of an abelian category.  
In particular, $\B$ has enough projectives and enough injectives and these coincide.  We denote the stable category by $\uB$, which is a triangulated category with the cosyzygy functor $\Omega^{-1}$ as its suspension \cite{TCRTA}.  In $\uB$ we will often write $X[i]$ for the $i^{th}$ suspension $\Omega^{-i}X$ of $X$.  We write $\ul{f}$ for the residue class in $\uB$ of a map $f$ in $\B$.  We further assume that all the Hom-spaces $\uB(X,Y)$ are finite-dimensional over $k$.  Typically, we have in mind for $\B$ either (an exact subcategory of) $\rmod A$ for a finite-dimensional self-injective $k$-algebra $A$ or else the category $\CM(R)$ of maximal Cohen-Macaulay modules over an isolated Gorenstein singularity $R$ (containing $k$).

For a subcategory $\C$ of $\B$, recall that a {\it right $\C$-approximation} of $X \in \B$ consists of a map $f: C_0 \rightarrow X$ with $C_0 \in \C$ such that any map $h : C \rightarrow X$ with $C \in \C$ can be factored through $f$.  The notion of  a {\it left $\C$-approximation} $g: X \rightarrow C_0$ is defined dually.  The subcategory $\C$ is said to be {\it functorially finite} in $\B$ if each object of $\B$ has both right and left $\C$-approximations.  Note that this condition is equivalent to requiring that the functors $\B(-,X)|_{\C}$ and $\B(X,-)|_{\C}$ are finitely generated (as functors from $\C$ to $\rmod k$) for each $X \in \B$.  Following Iyama \cite{Iyama1}, we say that a functorially finite subcategory $\C$ of $\B$ is  {\it maximal $(d-1)$-orthogonal} if 

\begin{eqnarray} \C & = & \{X \in \B\ |\ \uB(X,\C[i]) = 0, \forall\ 1\leq i < d\} = \{Y \in \B\ |\ \uB(\C,Y[i]) = 0, \forall\ 1\leq i < d\}. \end{eqnarray}  

  We shall henceforth assume that $\C$ is a functorially finite, maximal $(d-1)$-orthogonal subcategory of $\B$ for some $d \geq 1$.  In particular, $\C$ must contain all the projectives in $\B$ and we have $\uB(\C,\C[i]) = 0$ for all $1 \leq i < d$.  It is also easy to see that the induced subcategory $\uC$ of $\uB$ remains functorially finite and maximal orthogonal, and thus we may also view
$\uC$ as a maximal $(d-1)$-orthogonal subcategory of $\uB$.  If $\C = \add(T)$ for an object $T \in \B$, then we say that $T$ is a {\it $d$-cluster tilting object} (in $\B$ or in $\uB$).  Notice that in this case $\C$ will automatically be functorially finite.  Indeed, $\uC$ will be a finite type subcategory of $\uB$, of which we are assuming the Hom-spaces are finite-dimensional over $k$.  Thus, any $X \in \B$ has a right $\uC$-approximation $\ul{f} : C_0 \rightarrow X$ in $\uB$.  Then the map $(f\ p): C_0 \oplus P \rightarrow X$, 
where  $p : P \rightarrow X$ is a projective cover of $X$ in $\B$, gives a right $\C$-approximation of $X$.  The existence of left $\C$-approximations is established dually.

  We point out that for $d=1$ this definition forces $\C = \B$, which brings us back essentially to the setting considered by Auslander and Reiten in \cite{SEAA} and Buchweitz in \cite{Buch}.  With $\C$ and $d$ fixed we also define subcategories 
\begin{eqnarray}
\E_j & = & \{ X \in \B\ |\ \uB(\C,X[i]) = 0\ \mbox{for}\ 1 \leq i \leq d-1\ \mbox{and}\ i \neq j\}
\end{eqnarray}
for each $1 \leq j \leq d$.  Notice that $\E_{d} = \C$ and $\C \cup \C[1] \subseteq \E_{d-1}$.  If $d=2$,  then the defining condition for $\E_1$ becomes vacuous, and so in this case we set $\E_1 = \B$.

Our main results require an additional stronger vanishing condition on $\C$.  Fortunately, it turns out to be equivalent to a more natural (and more easily checked) periodicity condition, as we now verify.

\begin{lemma}  For $\C$ and $\B$ as above, the following are equivalent.
\begin{enumerate}
\item $\uB(\C,\C[i]) = 0$ for all $i$ with  $-d < i \leq -1$.
\item $\uC[d] = \uC$; that is, $\Omega^{d}C \in \C$ for each $C \in \C$.
\end{enumerate}
\end{lemma}

\noindent
{\it Proof.}  For $X \in \C$, notice that $X[d] \in \C$ if and only if $\uB(X[d],\C[i]) = 0$ for $1 \leq i < d$, which is equivalent to $\uB(X,\C[j]) = 0$ for $-d < j \leq -1$. $\Box$ \\

We will often assume that $\C$ satisfies the two equivalent conditions of the above lemma.  Note that these are automatic for $d=1$ and $\C = \B$.  In case $\uB$ has Serre duality $\uB(X,SY) \cong D\uB(Y,X)$ for an auto-equivalence $S$ of $\uB$, with $D$ denoting the duality $\Hom_k(-,k)$, then the above conditions are easily seen to be equivalent to $S(\C) = \C$.  

The following lemma is useful for obtaining exact sequences in $\B$, which may fail to be an abelian category.  It implies, in particular, that $\B$ has {\it plenty of projectives} in the terminology of \cite{Buch}.

\begin{lemma} For any map $f : X \rightarrow Y$ in $\B$, there exists an object $Z$ and a projective $P$ in $\B$ such that $\ses{Z}{X\oplus P}{Y}{\bnc{g}{i}}{(f\ p)}$ is exact in $\B$.  Moreover, there is a distinguished triangle $Z \stackrel{\ul{g}}{\rightarrow} X \stackrel{\ul{f}}{\rightarrow} Y \rightarrow$ in $\uB$, which determines $Z$ and $\ul{g}$ up to isomorphism in $\uB$.

\end{lemma}

\noindent
{\it Proof.}  Forming the pull-back of the exact sequence $\ses{\Omega Y}{P}{Y}{}{}$, where $P$ is projective, with respect to the map $f : X \rightarrow Y$ yields a commutative diagram in which the rows are exact sequences in $\B$:
$$\xymatrix{0  \ar[r] & \Omega Y \ar@{=}[d] \ar[r] & Z \ar[d] \ar[r] & X \ar[d]^f \ar[r] &0 \\ 0 \ar[r] & \Omega Y \ar[r] & P \ar[r]^p & Y \ar[r] & 0}.$$
Thus the sequence $\ses{Z}{X \oplus P}{Y}{}{(f\ p)}$ from the pull-back square is the one we want.  The second claim now follows from Lemma 2.7 in \cite{TCRTA} and the axioms for triangulated categories.  $\Box$\\


We use the standard notation $\rmod \C$ and $\rmod \uC$ for the categories of finitely presented contravariant $k$-linear functors from $\C$ and $\uC$, respectively, to $\rmod k$.   We also write $\stmod \uC$ for the stable category obtained from $\rmod \uC$ by factoring out the ideal of morphisms that factor through a projective.  As we only consider functors on $\C$ or $\uC$, and never on $\B$, all representable functors $\B(-,X)$ or $\uB(-,X)$ are to be interpreted as restricted to $\C$, and we forgo writing $\B(-,X)|_{\C}$ for the restriction.  We observe that our assumptions guarantee that all such representable functors  belong to $\rmod \C$ and $\rmod \uC$, respectively.  Indeed, we may complete a right $\uC$-approximation $\ul{f} : C_0 \rightarrow X$ to a triangle $Y \stackrel{\ul{g}}{\longrightarrow} C_0 \stackrel{\ul{f}}{\longrightarrow} X \rightarrow$, and then take a right $\uC$-approximation $\ul{h} : C_1 \rightarrow Y$.  This construction yields a projective presentation 
\begin{eqnarray} \uB(-,C_1) \stackrel{\uB(-,\ul{gh})}{\longrightarrow} \uB(-,C_0) \stackrel{\uB(-,\ul{f})}{\longrightarrow} \uB(-,X) \rightarrow 0. \end{eqnarray}
  Moreover, by the preceding lemma, these triangles may be lifted to short exact sequences $$\ses{Y \oplus Q}{C_0\oplus P_0}{X}{\bnc{g\ *}{*\ *}}{(f\ p)} \ \ \ \mbox{and}\ \ \ \ses{Z}{C_1 \oplus P_1}{Y \oplus Q}{}{\bnc{h\ *}{*\ *}}$$ with $P_0, P_1$ and $Q$ projective, which also yield right $\C$-approximations of $X$ and $Y \oplus Q$ resectively.  Splicing together the induced exact sequences of representable functors yields a projective presentation 
  \begin{eqnarray} \B(-,C_1\oplus P_1) \stackrel{\B(-,\varphi)}{\longrightarrow} \B(-,C_0 \oplus P_0) \stackrel{\B(-,(f\ p))}{\longrightarrow} \B(-,X) \rightarrow 0\end{eqnarray}
   where $\varphi$ has the form $\scriptsize \left( \begin{array}{cc} gh & * \\ * & * \end{array}\right)$.  Furthermore, we can now see that the representable functor $\uB(-,X)$ is also in $\rmod \C$ since it arises as the cokernel of the map $\B(-,P_X) \stackrel{\B(-,\pi_x)}{\longrightarrow} \B(-,X)$ induced by the projective cover $\pi_X : P_X \rightarrow X$.

Our current goal is to describe the projective resolutions of finitely presented $\uC$-modules in both $\rmod C$ and $\rmod \uC$.  We start with a simple but important observation that generalizes a theorem of Buan, Marsh and Reiten for $2$-cluster tilting objects in cluster categories \cite{CTA}  (see also Corollary 6.4 in \cite{IY}).  For a subcategory $\A$ of $\B$ we write $\gen{A}$ for the ideal of $\B$ generated by the identity morphisms of the objects of $\A$.

\begin{lemma} Let $\B$ and $\C$ be as above, and assume $d \geq 2$.
\begin{enumerate}
\item For any $M \in \rmod \uC$, we have $M \cong \uB(-,X)$ for some $X \in \E_{d-1}$ (without projective summands).
\item The functor $\eta: \uB \longrightarrow \rmod \uC$ given by $\eta(X) = \uB(-,X)$ is full and dense.  Moreover, the restriction of $\eta$ to $\ul{\E_{d-1}}$ induces a category equivalence $$\eta: \ul{\E_{d-1}}/\gen{\uC[1]} \stackrel{\approx}{\longrightarrow} \rmod \uC.$$
\end{enumerate}
 In particular, if $\uB$ has finite type, then so does $\rmod \uC$.
\end{lemma}

\noindent
{\it Proof.}  A minimal projective presentation of $M$ in $\rmod \uC$ has the form 
\begin{eqnarray}
\uB(-,C_1) \stackrel{\uB(-,f)}{\longrightarrow} \uB(-,C_0) \longrightarrow M \rightarrow 0 \end{eqnarray}
for a map $f : C_1 \rightarrow C_0$ in $\C$.  We can complete $\ul{f}$ to a triangle $C_1 \stackrel{\ul{f}}{\longrightarrow} C_0 \stackrel{\ul{g}}{\longrightarrow} X \longrightarrow$ in $\uB$.  The long-exact Hom-sequence now yields the exact sequence (using $d \geq 2$)
\begin{eqnarray} \uB(-,C_1) \stackrel{\uB(-,f)}{\longrightarrow} \uB(-,C_0) \stackrel{\uB(-,g)}{\longrightarrow} \uB(-,X) \longrightarrow \uB(-,C_1[1]) = 0,
\end{eqnarray}
 whence $M \cong \uB(-,X)$.  Furthermore, the exact sequences $$0=\uB(-,C_0[i]) \longrightarrow \uB(-,X[i]) \longrightarrow \uB(-,C_1[i+1])=0$$ for $1 \leq i \leq d-2$ show that $X \in \E_{d-1}$.

It follows easily that $\eta$ (even restricted to $\ul{\E_{d-1}}$) is full and dense, so we need only compute its kernel on $\ul{\E_{d-1}}$.  Clearly the kernel contains the ideal $\gen{\uC[1]}$ since $\uB(-,C[1]) = 0$ for all $C \in \C$.  Now let $f : X \rightarrow Y$ be a map between two objects of $\ul{\E_{d-1}}$ such that $\uB(C,f) = 0$ for all $C \in \C$.  If we complete a right $\uC$-approximation $g: C_0 \rightarrow X$ to a triangle $Z \longrightarrow C_0 \longrightarrow X \rightarrow$ in $\uB$, then the induced long exact sequence of representable functors on $\C$ shows that $Z \in \E_{d} = \C$.  As $fg = 0$ by assumption, we know that $f$ must factor through the connecting morphism $X \rightarrow Z[1]$, whence $f$ is in the ideal generated by $\uC[1]$.  $\Box$\\

\noindent
{\bf Remark.}  Of course, the final statement fails for $d=1$ as it is well-known that the stable Auslander algebra of a self-injective algebra of finite representation type usually has infinite representation type.\\

Before going on, we pause briefly to review some basics about finitely-presented functors and to explain some of our notation.  These facts are essentially due to Auslander and Reiten \cite{SEDRV}, but we shall follow the notation of \S 3 of \cite{Buch}.  Corresponding to the natural functor $p : \C \rightarrow \uC$, we have a restriction functor $p_* : \rmod \uC \rightarrow \rmod \C$, which is full and faithful and identifies $\rmod \uC$ with the full subcategory of $\rmod \C$ consisting of functors that vanish on projectives.  Moreover, $p_*$ has a right-exact left adjoint $p^*$ that is determined by $p^*\B(-,C) = \uB(-,C)$ for each $C \in \C$.  We interpret this functor, which takes $\C$-modules to $\uC$-modules, as tensoring with $\uC$ over $\C$, and we write $\Tor^{\C}_*(-,\uC)$ for its left derived functors.  Furthermore, by considering the projective presentations (2.4) and (2.3), we see that in fact $p^*\B(-,X) \cong \uB(-,X)$ for all $X \in \B$.

\begin{propos} Let $M \in \rmod \uC$, and assume that $d \geq 2$ and $\uC[d] = \uC$.
\begin{enumerate}
\item There is a projective presentation of $M$ in $\rmod \C$ of the form
$$\ses{\B(-,\Omega X) \longrightarrow \B(-,C_1)}{\B(-,C_0)}{M}{\B(-,f)}{}$$ for $C_0, C_1 \in \C$ and some $X \in \B$ with $M \cong \uB(-,X)$.  

\item Via $p^*$, the above sequence induces the following projective presentation of $M$ in $\rmod \uC$
$$\ses{\uB(-,\Omega X) \longrightarrow \uB(-,C_1)}{\uB(-,C_0)}{M}{\uB(-,f)}{}.$$
\item For any $X \in \ul{\E_{d-1}}$ we have a natural isomorphism  $\Omega^2_{\uC}[\uB(-,X)] \cong \uB(-,\Omega X)$ in $\stmod \uC$.

\end{enumerate}
\end{propos}

\noindent
{\it Proof.}  As in the preceding proof we can find $X \in \E_{d-1}$ with $M \cong \uB(-,X)$.  For simplicity, we assume that $X$ has no projective summands.  Keeping the notation introduced above and continuing the sequence (2.5) to the left, we obtain the exact sequence $$\ses{\uB(-,\Omega X) \longrightarrow \uB(-,C_1)}{\uB(-,C_0)}{M}{\uB(-,f)}{}$$ as $\uB(-,C_0[-1])=0$.  This sequence establishes (2) and also induces the isomorphism in (3), which can be seen to be natural in $X \in \ul{\E_{d-1}}$.  Using Lemma 2.2 we now lift the triangle $C_1 \stackrel{\ul{f}}{\longrightarrow} C_0 \stackrel{\ul{g}}{\longrightarrow} X \longrightarrow$ to a short exact sequence $\ses{C_1 \oplus P_1}{C_0 \oplus P_0}{X}{}{(g\ p)}$ in $\B$ with $P_0, P_1$ projective.  Notice that $(g\ p)$ is a right $\C$-approximation, since $\ul{g}$ is a right $\uC$-approximation by (2.6).  It follows that 
\begin{eqnarray} \ses{\B(-,C_1 \oplus P_1)}{\B(-,C_0 \oplus P_0)}{\B(-,X)}{}{} \end{eqnarray}
 is a projective resolution of $\B(-,X)$ in $\rmod \C$.  Taking a projective cover $\pi_X$ of $X$, the short exact sequence $\ses{\Omega X}{P_X}{X}{}{\pi_X}$ yields the exact sequence 
 \begin{eqnarray}
 \ses{\B(-,\Omega X) \longrightarrow \B(-,P_X)}{\B(-,X)}{\uB(-,X)}{\B(-,\pi_X)}{}
 \end{eqnarray} 
in $\rmod \C$.  Writing $\p(-,X)$ for the image of $\B(-,\pi_X)$, we can obtain the projective presentation of $M \cong \uB(-,X)$ as the mapping cone of the map from the sequence $$\ses{\B(-,\Omega X)}{\B(-,P_X)}{\p(-,X)}{}{}$$ to the sequence (2.7) which is induced by the inclusion $\p(-,X) \rightarrow \B(-,X)$.  Renaming $C_0 := C_0 \oplus P_0$ and $C_1 := C_1 \oplus P_1 \oplus P_X$ we see that this mapping cone has the desired form as in (1). $\Box$\\

\noindent
{\bf Remark.}  If $d=1$ and $\C = \B$, then the entire projective resolution of any $M = \uB(-,X)$ in $\rmod \C$ has the form (2.8) (cf. \cite{SEAA, SEDRV}), which is an instance of the presentation in part (1) of the proposition.  Thus, part (1) remains true in case $d=1$.  On the other hand, parts (2) and (3) of the proposition do not have interesting analogues in this case, since $M = \uB(-,X)$ will be projective in $\rmod \uC$.   Part (1), however, would yield a natural isomorphism $\Omega^2_{\C}[\uB(-,X)] \cong \B(-,\Omega X)$ in $\stmod \C$ for any $X \in \uB$, which resembles the isomorphism in (3).\\


We now describe the remaining terms of these projective resolutions for arbitrary $d \geq 2$.  Unfortunately, $X \in \E_{d-1}$ usually does not imply $\Omega X \in \E_{d-1}$, and hence we cannot simply repeat the above construction to build a projective resolution in $\rmod \uC$.  However, we'll see that we can iterate the construction, once the first $d+2$ terms of the resolution have been found.

\begin{therm} Let $\C$ be a maximal $(d-1)$-orthogonal subcategory of $\B$ with $\uC[d] = \uC$ and $d \geq 2$, and let $M \in \rmod \uC$.
\begin{enumerate}
\item  $M$ has a projective resolution in $\rmod \C$ of the form
$$\ses{\B(-,C_{d+1}) \rightarrow \cdots \rightarrow \B(-,C_1)}{\B(-,C_0)}{M}{}{}$$ with each $C_i \in \C$.
\item The induced sequence of functors on $\uC$ $$\ses{\Tor^{\C}_{d+1}(M,\uC)}{\uB(-,C_{d+1}) \rightarrow \cdots \rightarrow \uB(-,C_0)}{M}{}{}$$ is exact, and hence yields the first $d+2$ terms of a projective resolution for $M$ in $\rmod \uC$.
\item $\Tor_i^{\C}(M,\uC) = 0$ for all $i \neq 0, d+1$.
\item We have isomorphisms $\Tor^{\C}_{d+1}(M,\uC) \cong \Omega^{d+2}_{\uC} (M)$ in $\stmod \uC$ which are natural in $M$.
\item For any $X \in \ul{\E_{d-1}}$, we have a natural isomorphism $\Omega^{d+2}_{\uC}[\uB(-,X)] \cong \uB(-,\Omega^{d}X)$ in $\stmod \uC$.
\end{enumerate}
\end{therm}

\noindent
{\it Proof.}  As in Proposition 2.3, there is a triangle $C_1 \rightarrow C_0 \rightarrow X \rightarrow$ in $\uB$ with $M \cong \uB(-,X)$ and $X \in \E_{d-1}$.  Thus $\Omega X = X[-1] \in \E_1$.  We set $L_1 : = \Omega X$, and recursively define $L_j$ for $j \geq 2$ as follows: Take a right $\uC$-approximation $f_j : C_j \rightarrow L_{j-1}$ and complete it to a triangle $L_j \longrightarrow C_j \stackrel{f_j}{\longrightarrow} L_{j-1} \rightarrow$ in $\uB$.

  We prove by induction that 
  \begin{itemize}
  \item[(i)] $L_j \in \E_j$ for each $1 \leq j \leq d$; and 
  \item[(ii)] $\uB(-,L_j[j-d]) \cong \uB(-,X[-d])$ for $1 \leq j \leq d-1$.
  \end{itemize}
For $j=1$, we have already noted that (i) holds, and (ii) is trivial.  Now assume that both statements hold for some $j$ with $1 \leq j < d$.   We consider the exact sequences in $\rmod \uC$ for various $i$
$$\uB(-,L_j[i-1]) \longrightarrow \uB(-,L_{j+1}[i]) \longrightarrow \uB(-,C_{j+1}[i]).$$
By hypothesis, the first term vanishes for all $i$ with $2 \leq i \leq d$ and $i \neq j+1$; while the third term vanishes for all $i$ with $1 \leq i \leq d-1$.  We thus see that the middle term vanishes for all $i \neq j+1$ with $2 \leq i \leq d-1$.  It vanishes for $i=1$ since $f_j$ is a right $\uC$-approximation, making $\uB(-,f_j)$ surjective.  This establishes $L_{j+1} \in \E_{j+1}$.  In particular, observe that $C_{d+1} := L_{d} \in \E_{d} = \C$.  To see (ii), assume $j < d-1$ and notice that $\uB(-,L_{j+1}[j+1-d]) \cong \uB(-,L_j[j-d]) \cong \uB(-,X[-d])$ since $\uB(\C,C_{j+1}[i]) = 0$ for $i = j+1-d, j-d$.

For each $j$ with $1 \leq j \leq d-2$ we now have a short exact sequence
\begin{eqnarray} \ses{\uB(-,L_{j+1})}{\uB(-,C_{j+1})}{\uB(-,L_{j})}{}{}\end{eqnarray} in $\rmod \uC$ since $\uB(\C,L_j[-1]) \cong \uB(\C[d],L_j[d-1]) \cong \uB(\C,L_j[d-1])=0$ and $f_{j+1}$ is a right $\uC$-approximation.  Similarly, for $j = d-1$, the triangle $C_{d+1} \longrightarrow C_{d} \longrightarrow L_{d-1} \rightarrow$ induces an exact sequence $$\ses{\uB(-,L_{d-1}[-1]) \longrightarrow \uB(-,C_{d+1})}{\uB(-,C_d)}{\uB(-,L_{d-1})}{}{}.$$   Splicing these sequences together and using the isomorphism $\uB(-,L_{d-1}[-1]) \cong \uB(-,\Omega^{d}X)$ from (ii) yields an exact sequence 
\begin{eqnarray}
0 \rightarrow \uB(-,\Omega^d X) \rightarrow \uB(-,C_{d+1}) \rightarrow \cdots \rightarrow \uB(-,C_2) \rightarrow \uB(-,\Omega X) \rightarrow 0
\end{eqnarray}
in $\rmod \uC$, which can be viewed as the beginning of a projective resolution for $\uB(-,\Omega X)$.  Now splicing (2.10) with the projective presentation from Proposition 2.4(2) gives the first $d+2$ terms of a projective resolution for $M$ in $\rmod \uC$.  The isomorphism in (5) follows, and its naturality is a routine verification.

  At the same time, applying Lemma 2.2 to each triangle $L_j \longrightarrow C_j \stackrel{f_j}{\longrightarrow} L_{j-1} \rightarrow$ we obtain exact sequences $\ses{L_j}{C_j \oplus P_j}{L_{j-1}}{}{}$ in $\B$ and exact sequences $\ses{\B(-,L_j)}{\B(-,C_j \oplus P_j)}{\B(-,L_{j-1})}{}{}$ in $\rmod \C$.  Splicing these together, we obtain a projective resolution for $\B(-,\Omega X)$ in $\rmod \C$
\begin{eqnarray}
\ses{\B(-,C_{d+1})}{\B(-,C_{d} \oplus P_{d}) \rightarrow \cdots \rightarrow \B(-,C_2 \oplus P_2)}{\B(-,\Omega X)}{}{}.
\end{eqnarray}
Combining this with the projective presentation in Proposition 2.4, yields the desired resolution of $M$.  If we now apply $-\otimes_{\C} \uC$ to this resolution, the exactness of (2.10) and of $\ses{\uB(-,\Omega X)}{\uB(-,C_1) \longrightarrow \uB(-,C_0)}{M}{}{}$ shows that $\Tor_i^{\C}(M,\uC) = 0$ for all $i \neq 0, d+1$, and $\Tor_{d+1}^{\C}(M,\uC) \cong \Omega^{d+2}_{\uC}(M)$.  Moreover, this last isomorphism is clearly natural in $M$.  $\Box$\\

If $M = \uB(-,C)$ for a nonprojective $C \in \C$, then the projective resolution in $\rmod \C$ from the above theorem takes on an even simpler form.  As in Proposition 2.4, the second syzygy of $M$ is isomorphic to $\B(-,\Omega C)$.  Since $\uB(\C,\Omega C) = 0$ the projective cover of $\Omega C$ will be a right $\C$-approximation.  We thus obtain an exact sequence $\ses{\B(-,\Omega^2C)}{\B(-,P_2)}{\B(-,\Omega C)}{}{}$ in $\rmod \C$ with $P_2$ projective.  Repeating this construction, using $\uB(\C,\Omega^i C)=0$ for $1 \leq i \leq d-1$, we obtain the projective resolution:
$$\ses{\B(-,\Omega^{d}C)}{\B(-,P_{d}) \rightarrow \cdots \rightarrow \B(-,P_2)}{\B(-,P_C) \longrightarrow \B(-,C) \longrightarrow \uB(-,C)}{}{}$$
with $\Omega^{d}C \in \C$ by assumption.  Passing to $\uB$ by factoring out the maps that factor through projectives, all terms of this projective resolution vanish except for the $0^{th}$ and $(d+1)^{th}$ terms.  In particular, we recover the following isomorphisms 
\begin{eqnarray} \Tor_{d+1}^{\C}(\uB(-,C),\uC) \cong \uB(-,\Omega^{d}C)
\end{eqnarray}
of functors on $\uC$, which are natural in $C \in \uC$ (note that they also follow from combining parts (4) and (5)).  Thus we have isomorphisms of bifunctors on $\uC$
  \begin{eqnarray} \Tor_{d+1}^{\C}(\uB(-,-),\uC) \cong \uB(-,\Omega^{d}(-)).
\end{eqnarray}

We also point out that the remaining terms of the projective resolution of $\uB(-,X)$ in $\rmod \uC$ can now be obtained by essentially shifting the terms described in part (2) of the theorem, and in this way we obtain a {\it quasi-periodic} projective resolution for $\uB(-,X)$.  This is due to the assumption that $\uC[d] = \uC$, which guarantees that $\E_{i}[d] = \E_{i}$ for each $i$.  Hence $\Omega^{d+2}_{\C}[\uB(-,X)] \cong \uB(-,\Omega^d X)$ with $\Omega^d X = X[-d] \in \E_{d-1}$.  Then the construction from the proof can clearly be shifted by the $-d^{th}$ power of the suspension functor to obtain the next $d+2$ terms of the projective resolution:
$$\ses{\uB(-,X[-2d])}{\uB(-,C_{d+1}[-d]) \rightarrow \cdots \rightarrow \uB(-,C_0[-d])}{\uB(-,X[-d])}{}{},$$ and so on.  We also easily see that iterating the isomorphism from part (5) of the theorem yields isomorphisms $\Omega_{\uC}^{s(d+2)}[\uB(-,X)] \cong \uB(-,\Omega^{sd}X)$ in $\stmod \uC$ for each $s \geq 1$, which are natural in $X \in \ul{\E_{d-1}}$.

\section{Bimodule resolutions of stable Auslander algebras}
\setcounter{equation}{0}

In this section we specialize to the case where $\C = \add(T)$ for a $d$-cluster tilting object $T \in \B$ with $d \geq 1$.  The evaluation functor $ev_T : M \mapsto M(T)$ gives category equivalences $\rmod \C \rightarrow \rmod \Lambda$ and $\rmod \uC \rightarrow \rmod \Gamma$, where $\Lambda = \Endo_{\B}(T)$ and $\Gamma = \stEndo_{\B}(T)$.  Our $\Hom$-finiteness assumption on $\uB$ guarantees that $\Gamma$ is finite-dimensional, although $\Lambda$ need not be.  We also note that $\Gamma$ may be decomposable as an algebra, and may even have semisimple blocks which we typically want to ignore.  As we deal with bimodules, we assume for convenience that $k$ is perfect (although, it suffices to know that $\Gamma$ splits over a separable extension of $k$).  Under this assumption, the projective bimodule summands of $\Gamma$ correspond precisely to semisimple blocks.

  We now translate some of our above results (parts (3) and (4) of Theorem 2.5 and (2.13)) to this setting in the corollary below.  These statements are also true for $d=1$ by Theorem 1.1 and Proposition 6.5 of \cite{Buch}.

\begin{coro} Let $T \in \B$ be a $d$-cluster tilting object with $d \geq 1$ such that $\Omega^{d}T \cong T$ in $\uB$, and set $\Lambda = \Endo_{\B}(T)$ and $\Gamma = \stEndo_{\B}(T)$.  Then 
\begin{enumerate}
\item $\Tor_i^{\Lambda}(-,\Gamma) = 0$ on $\rmod \Gamma$ for all $i \neq 0, d+1$; 
\item $\Tor_{d+1}^{\Lambda}(-,\Gamma) \cong \Omega^{d+2}$ as functors on $\stmod \Gamma$.
\item $\Tor_{d+1}^{\Lambda}(\Gamma,\Gamma) \cong \uB(T,\Omega^{d}T)$ as $(\Gamma,\Gamma)$-bimodules.
\end{enumerate}
\end{coro}

The assumption that $\Omega^{d}T \cong T$ implies that $\uB(T,\Omega^{d}T)$ is isomorphic to a twisted bimodule ${}_{\sigma}\Gamma_1$ for some $k$-algebra automorphism $\sigma$ of $\Gamma$, which corresponds to an isomorphism $\eta: \Omega^{d}T \stackrel{\cong}{\longrightarrow} T$.  If $\Omega^{d} \cong Id$ as functors on $\add(T)$, then $\uB(T, \Omega^{d}T) \cong \Gamma$ as bimodules.

We now delve deeper to obtain information about the projective resolution of $\Gamma$ over its enveloping algebra $\Gamma^e$.  Recall that $\Gamma$ is {\it periodic} if this resolution is periodic.  We will also say that $\Gamma$ is {\it quasi-periodic} (or, equivalently, that this resolution is quasi-periodic) if $\Omega^n_{\Gamma^e}(\Gamma)$ is isomorphic to a twisted bimodule ${}_{\sigma}\Gamma_1$ as above.  In this case, it easily follows that each finitely generated $\Gamma$-module has bounded Betti numbers.

\begin{therm} Let $T \in \B$ be a $d$-cluster tilting object such that $\Omega^{d}T \cong T$ in $\uB$, and set $\Lambda = \Endo_{\B}(T)$ and $\Gamma = \stEndo_{\B}(T)$.  Then
\begin{enumerate}
\item $\Tor_{d+1}^{\Lambda}(-,\Gamma) \cong -\otimes_{\Gamma} \Tor_{d+1}^{\Lambda}(\Gamma,\Gamma)$ as functors on $\rmod \Gamma$.
\item $\Omega^{d+2}_{\Gamma^e}(\Gamma) \cong \Tor_{d+1}^{\Lambda}(\Gamma, \Gamma) \cong \uB(T,\Omega^{d}T)$ as $(\Gamma,\Gamma)$-bimodules (up to projective summands).
\end{enumerate} 
In particular, $\Gamma$ is self-injective.  Moreover, writing $\Gamma = \Gamma_0 \times \Gamma_{s}$ where $\Gamma_s$ is the largest semisimple direct factor of $\Gamma$, we see that $\Gamma_0$ is quasi-periodic of quasi-period $d+2$.  If $\Omega^{dr}|_{\add(T)} \cong Id_{\add(T)}$ as  functors for some $r \geq 1$, then $\Gamma_0$ is periodic with period dividing $r(d+2)$.
\end{therm}

\noindent
{\bf Remarks.} (1) Part (2) and its consequences can be viewed as an extension of Theorem 1.5 in \cite{Buch}.  Notice that we can avoid assuming that $\Lambda$ has Hochschild dimension $d+1$, even when $d=1$, since our broader assumptions on $\B$ and $T$ guarantee that $\Gamma$ is finite-dimensional and self-injective, and we will see that these conditions suffice.  In particular, this simplifies certain issues arising in applications of Buchweitz's results (cf. 1.6, 1.12 in \cite{Buch}).

(2) While quasi-periodicity appears weaker than periodicity, we are unaware of any finite-dimensional algebras that are quasi-periodic but not periodic.  This theorem could potentially be used to produce such examples: for instance, one would need a $d$-cluster tilting object $T$ with $\Omega^d T \cong T$ but where no positive power of $\Omega^d$ is isomorphic to the identity functor on $\add(T)$.\\

\noindent
{\it Proof.}  For (1), notice that $\Tor_{d+1}^{\Lambda}(-,\Gamma)$ is an exact functor on $\rmod \Gamma$ as $\Tor_{d}^{\Lambda}(-,\Gamma)=\Tor_{d+2}^{\Lambda}(-,\Gamma) = 0$.  Thus $\Tor_{d+1}^{\Lambda}(-,\Gamma) \cong -\otimes_{\Gamma} \Tor_{d+1}^{\Lambda}(\Gamma,\Gamma)$ by the Eilenberg-Watts theorem.  Observe that $\Tor_{d+1}^{\Lambda}(\Gamma,\Gamma) \cong {}_{\sigma}\Gamma_1$ is a projective $\Gamma$-module on either side.  Furthermore, since we have an invertible bimodule $\Tor_{d+1}^{\Lambda}(\Gamma,\Gamma)$ inducing $\Omega^{d+2}$ on $\stmod \Gamma$, we see that $\Omega$ must be an equivalence and $\Gamma$ is self-injective.

For (2), let $\cdots \rightarrow P_1 \stackrel{f_1}{\longrightarrow} P_0 \longrightarrow \Lambda \rightarrow 0$ be a projective resolution of $\Lambda$ over $\Lambda^e$.  Applying $-\otimes_{\Lambda^e} \Gamma^e$ yields a complex $Q_{\bullet} := \Gamma \otimes_{\Lambda} P_{\bullet} \otimes_{\Lambda} \Gamma$ of projective $\Gamma^e$-modules with homology given by $$\Tor_*^{\Lambda^e}(\Lambda, \Gamma^e) \cong \Tor_*^{\Lambda}(\Gamma, \Gamma).$$  As Corollary 3.1 tells us that this homology vanishes in all degrees except $0$ and $d+1$, the beginning of a projective resolution of $\Gamma$ over $\Gamma^e$ has the form $$\ses{\Omega^{d+2}(\Gamma) \oplus Q}{Q_{d+1} \rightarrow \cdots \rightarrow Q_0}{\Gamma}{}{}$$ for some projective bimodule $Q$.  Furthermore, from the definition of $\Tor$ we have an epimorphism\footnote{It is an isomorphism if $\Lambda$ has Hochschild dimension $d+1$.  This holds for instance if $\B = \rmod A$ for a finite-dimensional self-injective algebra $A$, as then $\Lambda$ is a finite-dimensional algebra of global dimension $d+1$ \cite{Iyama1}.} $\Omega^{d+2}(\Gamma) \oplus Q \rightarrow \Tor_{d+1}^{\Lambda}(\Gamma,\Gamma)$.  Let $K$ be the kernel and observe that $K$ is projective on either side since $\Tor_{d+1}^{\Lambda}(\Gamma,\Gamma)$ and $\Omega^{d+2}(\Gamma) \oplus Q$ both are.  Also observe that by definition $K = \im (1 \otimes f_{d+2} \otimes 1)$ consists of the $(d+1)$-boundaries of $Q_{\bullet}$.  We claim that $K$ is a projective $(\Gamma,\Gamma)$-bimodule; since $\Gamma$ is self-injective it will then follow that the short exact sequence $\ses{K}{\Omega^{d+2}(\Gamma) \oplus Q}{\Tor_{d+1}^{\Lambda}(\Gamma,\Gamma)}{}{}$ splits, yielding $\Omega^{d+2}(\Gamma) \cong \Tor_{d+1}^{\Lambda}(\Gamma,\Gamma)$ as bimodules (up to projective summands).


To see that $K$ is projective, we go back a step and apply $ \Gamma \otimes_{\Lambda} -$ to $P_{\bullet}$ to get a projective $(\Gamma,\Lambda)$-bimodule resolution $\Gamma \otimes _{\Lambda} P_{\bullet}$ of ${}_{\Gamma}\Gamma \otimes_{\Lambda} \Lambda_{\Lambda} \cong {}_{\Gamma} \Gamma_{\Lambda}$.  Set $L = \ker (1 \otimes f_{d+1}) \cong \coker (1 \otimes f_{d+3})$.  Since $-\otimes_{\Lambda} \Gamma$ is right-exact, we have $L \otimes_{\Lambda} \Gamma \cong \coker (1 \otimes f_{d+3} \otimes 1) \cong \im (1 \otimes f_{d+2} \otimes 1) = K$ as $(\Gamma,\Gamma)$-bimodules.  For any finitely-presented right $\Gamma$-module $M$, $M \otimes_{\Gamma} \Gamma \otimes_{\Lambda} P_{\bullet} \cong M \otimes_{\Lambda} P_{\bullet}$ is a projective resolution of $M_{\Lambda}$.  Since $\pdim\ M_{\Lambda} \leq d+1$, $M \otimes_{\Gamma} L \cong \coker (1_M \otimes f_{d+3}) \cong \ker (1_M \otimes f_{d+1})$ is a projective right $\Lambda$-module.  In particular, $M \otimes_{\Gamma} K \cong M \otimes_{\Gamma} (L \otimes_{\Lambda} \Gamma) \cong (M \otimes_{\Gamma} L) \otimes_{\Lambda} \Gamma$ is a projective right $\Gamma$-module for any $M$.  Since $K$ is projective on either side, Theorem 3.1 of \cite{TEG} implies that $K$ is a projective bimodule.

For the final statement, we may assume that $\Gamma$ has no semisimple blocks by working with $\Gamma_0$ and an appropriate summand $T_0$ of $T$ instead.  Observe that for any $r \geq 1$, $\Omega^{r(d+2)}(\Gamma) \cong \Omega^{d+2}(\Gamma)^{\otimes r} \cong \uB(T,\Omega^{d}T)^{\otimes r}$ up to projective summands by (2) and Corollary 3.1(3).  
Using part (1), Corollary 3.1(3) and (2.10) we now obtain $\uB(T,\Omega^{d}T)^{\otimes r} \cong \uB(T,\Omega^{rd}T)$  by induction on $r \geq 1$ (cf. Prop. 6.5 in \cite{Buch}).  Furthermore, the latter bimodule is isomorphic to $\Gamma = \uB(T,T)$ as a bimodule if and only if $\Omega^{rd}$ is isomorphic to the identity functor on $\add(T)$.  $\Box$\\

Many examples of cluster-tilting objects appear inside Calabi-Yau triangulated categories, such as the cluster categories of \cite{BMRRT} or categories of the form $\ul{\CM}(R)$ for an isolated Gorenstein hypersurface singularity $R$ \cite{BIKR}.  Recall that an auto-equivalence $S$ of $\uB$ is called a {\it Serre functor} if there exist natural isomorphisms $D\uB(X,Y) \cong  \uB(Y,SX)$ for all $X, Y \in \uB$, where $D = \Hom_k(-,k)$ is the duality with respect to the ground field.  In this case, there is a canonical enhancement of $S$ into a triangulated functor, and if $S \cong -[s]$ as triangulated functors on $\uB$, then we say that $\uB$ is {\it Calabi-Yau of dimension $s$}.  Here we will consider the weaker requirement that $S \cong -[s]$ only as $k$-linear functors, in which case we say that $\uB$ is {\it weakly Calabi-Yau of dimension $s$}, in the sense of \cite{CYTC}.  This amounts to the existence of natural isomorphisms $$D \uB(X,Y) \cong \uB(Y,X[s])$$ for all $X, Y \in \uB$ (In order for $\uB$ to be Calabi-Yau of dimension $s$, one additionally requires that these natural isomorphisms are compatible with the suspension functor as in Proposition 2.2 of \cite{CYTC}). 

In case $\uB$ is weakly $s$-Calabi-Yau, the injective objects in $\rmod \uB$ have the form $D\uB(X,-) \cong \uB(-,X[s]) \cong \uB(-[-s],X)$ for $X \in \uB$, which shows that $\rmod \uB$ is a Frobenius category with Nakayama equivalence given by $\nu : F \mapsto F \circ [-s]$.  Thus $\stmod \uB$ is a Hom-finite triangulated category.  Moreover, Serre duality in $\uB$ guarantees that $\uB$ is a dualizing $k$-variety in the sense of \cite{SEDRV}, and hence the Auslander-Reiten formula implies 
$$D\stHom_{\uB}(F,G) \cong \Ext^1_{\uB}(G, D\Tr F) \cong \stHom_{\uB}(G, \Omega_{\uB} \nu F)$$  for all $F, G \in \stmod \uB$; that is, $\Omega_{\uB} \nu : F \mapsto \Omega_{\uB}( F \circ [-s])$ is a Serre functor for $\stmod \uB$.  Moreover, knowledge of the projective resolution for $F \in \rmod \uB$ (from \cite{SEDRV} or \cite{SEAA}, for example) implies that $\Omega_{\uB}^3(F) \cong F \circ [1]$.  Hence $\nu \cong \Omega_{\uB}^{-3s}$ on $\stmod \uB$, and the Serre functor for $\stmod \uB$ satisfies $S = \Omega_{\uB} \nu \cong \Omega_{\uB}^{-(3s-1)}$, showing that $\stmod \uB$ is weakly $(3s-1)$-Calabi-Yau when $\uB$ is weakly $s$-Calabi-Yau (this has been observed elsewhere: see \cite{TOC}, for instance).  This result can in fact be viewed as the $d=1$ case of the following more general statement regarding maximal $(d-1)$-orthogonal subcategories of Calabi-Yau triangulated categories.  In the second part, we apply Theorem 3.2 to obtain a partial generalization of Proposition 2.1 in \cite{SCY}.

\begin{propos}[Cf. 5.4 in \cite{GKO}] Let $\C$ be a maximal $(d-1)$-orthogonal subcategory of $\B$ with $\uC[d] = \uC$, and assume that $\uB$ is weakly $sd$-Calabi-Yau for some integer $s$.  
\begin{enumerate}
\item $\stmod \uC$ is a weakly Calabi-Yau triangulated category of dimension $s(d+2)-1$.
\item If $\C = \add(T)$ for a $d$-cluster tilting object $T \in \B$ and $\Gamma = \stEndo_{\B}(T)$ has no semisimple blocks, then $\Omega^{-s(d+2)}_{\Gamma^e}(\Gamma) \cong D\Gamma$ as bimodules.
\end{enumerate}
\end{propos}

\noindent
{\it Proof.} (1) As remarked after Lemma 2.1, the assumption $\uC[d] = \uC$ ensures that $\uC$ is invariant under the Serre functor $S$ of $\uB$.  Hence the same argument given above for $\uB$ shows that $\rmod \uC$ is a Frobenius category with Nakayama equivalence $\nu$ given by $F \mapsto F \circ [-sd]$.  If $F = \uB(-,X) \in \rmod \uC$ for $X \in \E_{d-1}$, then $\nu(F) \cong \uB(-,X[sd]) \cong \Omega_{\uC}^{-s(d+2)}(F)$ by Theorem 2.5(5).  Since $\uC$ is also a dualizing $k$-variety (one again uses the Serre duality to check that the duality $D$ preserves finitely presented functors on $\uC$ and $\uC^{\op}$), the above argument also shows that a Serre functor for $\stmod \uC$ is given by $S = \Omega_{\uC} \nu \cong \Omega_{\uC}^{1-s(d+2)}$, and the claim follows.

(2) By Theorem 3.2, we have $\Omega^{-s(d+2)}(\Gamma) \cong \uB(T,\Omega^{-sd}T) \cong \uB(T,T[sd]) \cong D\uB(T,T) \cong D\Gamma$ as bimodules. $\Box$\\

\noindent
{\bf Remarks.}  (1) We point out that the curious requirement that the weak Calabi-Yau dimension of $\uB$ is $sd$ does not impose an unnecessary restriction in light of the assumption $\uC[d] = \uC$.  Indeed, if $\uB$ is weakly $n$-C.Y. then $\uB(C,C[n]) \cong D\uB(C,C) \neq 0$ for any $C \in \C$ implies that $d \mid n$.

(2) In fact, the full Calabi-Yau property is shown to hold for $\stmod \uC$ in \S 5 of \cite{GKO}, since $\uC$ with suspension $-[d]$ is a $(d+2)$-angulated category.

\section{Examples and concluding remarks}
\setcounter{equation}{0}

As remarked in the introduction, this work is motivated by the recent discovery of symmetric algebras with $D\Tr$-periodic module categories arising as stable endomorphism rings of $2$-cluster tilting objects in the Cohen-Macaulay module categories of $1$-dimensional hypersurface singularities \cite{BIKR}.  We briefly recall the construction introduced there, as we now know that it provides a powerful tool for producing periodic symmetric algebras of period $4$.  

Set $S = k[[x,y]]$ and $\mathfrak{m} = (x,y)$.  Choose irreducible power series $f_i \in \mathfrak{m} \setminus \mathfrak{m}^2$ for $1 \leq i \leq n$ with $(f_i) \neq (f_j)$ for $i \neq j$, and set $f = f_1 f_2 \cdots f_n$.  Then $R = S/(f)$ is an isolated hypersurface singularity of dimension $1$, and $T = \oplus_{i=1}^n S/(f_1 \cdots f_i)$ is a $2$-cluster tilting object in $\CM (R)$.  Moreover, Eisenbud's matrix factorization theorem implies that $\Omega^2 \cong Id$ on $\sCM(R)$, and thus on $\add(T)$ as well.  Hence Theorem 3.2 implies that $\Gamma = \stEndo_R(T)$ is periodic of period $4$.  The quiver of $\Gamma$ (but not the relations) is described in Proposition 4.10 of \cite{BIKR}:
$$\xymatrix{1  \ar[r]<0.5ex> & 2  \ar[r]<0.5ex> \ar[l]<0.5ex> & \cdots \ar[l]<0.5ex> \ar[r]<0.5ex> & n-2 \ar[l]<0.5ex> \ar[r]<0.5ex> & n-1 \ar[l]<0.5ex>}$$
with a loop at vertex $i$ if and only if $(f_i, f_{i+1}) \neq \mathfrak{m}$.  Furthermore, it is shown that two families of algebras of quaternion type  are explicitly realized in this way.  These algebras are known to have tame representation type, but starting with a hypersurface $R$ of wild CM-type should produce an algebra $\Gamma$ of wild type and period $4$.

Our results also yield new information in the classical case where $d=1$.  For example, if $R$ is a simple curve singularity of finite CM-type (in arbitrary characteristic) and $\Gamma$ is the stable Auslander algebra of $\CM (R)$, it follows from Theorem 1.1(3) that $\Gamma$ is periodic of period dividing $6$.  Moreover, since $\sCM(R)$ is $2$-Calabi-Yau, $\stmod \Gamma$ will be (weakly) $5$-Calabi-Yau by Proposition 3.3.  The algebras $\Gamma$ that arise in this way are (a proper subset of the) deformed preprojective algebras of generalized Dynkin type, as studied in \cite{ErdSko2}.  We have previously applied this information about the periods and stable Calabi-Yau dimensions of these algebras in the study of the same properties for the representation-finite self-injective algebras \cite{SCY}.  Similarly, if $R$ is a two-dimensional simple surface singularity (in arbitrary characteristic), the stable Auslander algebra $\Gamma$ of $\CM (R)$ is periodic of period dividing $6$ and stably $2$-Calabi-Yau.  The algebras $\Gamma$ arising in this way are necessarily deformed preprojective algebras of Dynkin type by \cite{BES}, and it is an interesting problem whether every such deformed preprojective algebra is isomorphic to the stable Auslander algebra of $\CM (R)$ for some simple surface singularity $R$ in arbitrary characteristic, as classified in \cite{GK}.

\vspace{2mm}
Unfortunately, it is still a challenging problem to find additional examples of maximal $(d-1)$-orthogonal subcategories where our results can be applied.  For instance, Erdmann and Holm \cite{ErdHolm} have shown that maximal $(d-1)$-orthogonal subcategories rarely exist in $\B = \rmod A$ for a self-injective $k$-algebra $A$.  Specifically, they show that they can only exist if every finite-dimensional $A$-module has complexity at most $1$.  Such algebras do exist -- periodic algebras, for example -- but even here the examples are limited.  Known examples of periodic algebras include all self-injective algebra of finite representation type \cite{Per}, but any periodic algebra constructed as the stable endomorphism ring of a maximal $(d-1)$-orthogonal subcategory in this context, will again have finite representation type by Lemma 2.3.  Still, it would be interesting to see which self-injective algebras of finite representation type are $d$-cluster tilted in this sense.   One could also look for maximal $(d-1)$-orthogonal subcategories of modules over tame and wild periodic algebras, which include the algebras of quaternion type, the preprojective algebras of Dynkin type and the $m$-fold mesh algebras \cite{ErdSko2}.  

Nevertheless, it may still be possible to find interesting examples of $d$-cluster tilting objects in {\it subcategories} of stable module categories.  In particular, our main results can be applied to a (finite type) maximal $(d-1)$-orthogonal subcategory inside some exact Frobenius subcategory $\B$ of $\rmod A$.  Namely, in light of Erdmann and Holm's result, one should take $\B$ to be the full subcategory of $\rmod A$ consisting of modules of complexity at most $1$, which is an exact subcategory with $\uB$ a triangulated subcategory of $\stmod A$.  Even here, however, it is not clear whether one will be able to find a module satisfying the restrictive self-orthogonality and Ext-configuration conditions required of a cluster-tilting object. 

\vspace{2mm}
Another source of applications can be found in the exciting work of Iyama and Oppermann on {\it higher preprojective algebras} \cite{IO}.  If $A$ is a finite-dimensional algebra with $\gldim A \leq n$ for which $\rmod A$ contains an $n$-cluster-tilting object, then the $(n+1)$-preprojective algebra of $A$ can be defined as $\tilde{A} = T_A \Ext^n_A(DA,A)$, the tensor algebra over $A$ of the bimodule $\Ext^n_A(DA,A)$.  Moreover, Iyama and Oppermann show  that $\tilde{A}$ can be realized as the endomorphism ring of an $n$-periodic $n$-cluster-tilting object in a certain Hom-finite triangulated category (namely, the $n$-Amiot cluster category $\C_{A}^n$ associated to $A$).  It follows immediately from Theorem 3.2 that $\tilde{A}$ has at least a quasi-periodic projective resolution over its enveloping algebra.  However, it appears a nontrivial problem to determine the order of the $n^{th}$ shift functor $[n]$ on the relevant maximal $(n-1)$-orthogonal subcategory of $\C_A^n$, and thus to determine whether or not this resolution is indeed periodic.

  For example, if $n=1$ and $A$ is a hereditary algebra of finite representation type, then the corresponding $2$-preprojective algebra will be the usual preprojective algebra associated to the (Dynkin) quiver of $A$.  Here $\tilde{A}$ is the endomorphism ring of a $1$-periodic $1$-cluster tilting object $T$, but has period $6$ (with some exceptions in characteristic $2$ where the period is $3$).  This means that for the $T$ in question, one has $T[1] \cong T$ but $-[1] : \add(T) \rightarrow \add(T)$ is not isomorphic to the identity functor, although its square $-[2]$ is.
  
  A more interesting example with $n=2$ can be found in \cite{IO}, Example 4.18.   Here we have a $3$-preprojective algebra $\tilde{A}$ for which $\Omega^{12}$ fixes each simple module up to isomorphism.   Since $\tilde{A}$ is the endomorphism ring of a $2$-periodic $2$-cluster-tilting object $T$, with $-[2]|_{\add(T)}$ inducing $\Omega^4$ on $\stmod \tilde{A}$, we see that the order of $-[2]$ on $\add(T)$ must be a multiple of $3$ (if it is finite).    
  
  Finally, we point out that Proposition 3.3 applies to all of the $(n+1)$-preprojective algebras $\tilde{A}$, since the relevant $n$-Amiot cluster category is $n$-Calabi-Yau by construction.  Thus part (2) of the proposition shows that $\Omega^{-n-2}_{\tilde{A}^e}(\tilde{A}) \cong D\tilde{A}$ as bimodules.  Since $D\tilde{A} \cong {}_1 \tilde{A}_{\nu}$ for the Nakayama automorphism $\nu$ of $\tilde{A}$, we can see that $\tilde{A}$ is periodic if and only if $\nu$ has finite order in the group of outer automorphisms of $\tilde{A}$.  However, even this latter condition remains difficult to verify.

\end{document}